\documentclass[11pt]{amsart}
\usepackage{amsmath}
\usepackage{amssymb}
\usepackage{mathpazo}
\usepackage{eucal}
\usepackage{amsthm}
\usepackage{amscd}
\usepackage{hyperref}
\usepackage{url}
\usepackage{skull}
\usepackage{adjustbox}

\usepackage{tikz}
\usetikzlibrary{shapes.geometric, arrows}
\usetikzlibrary{shapes,arrows}
\usepackage[latin1]{inputenc}
\usepackage{tikz}
\usetikzlibrary{shapes,arrows}
\usetikzlibrary{shapes.multipart}
\usepackage{verbatim}
\usepackage{tkz-berge}
\usepackage{tikz-qtree}

\usepackage{amssymb,}
\usepackage[]{amsmath, amsthm, amsfonts,graphicx, amscd,}
\usepackage[all,cmtip]{xy}
\input amssym.def \input amssym
\usepackage{forest}
\usepackage{mdframed}
\usepackage{framed}

\begin{document}

\newtheorem{thm}{Theorem}[section]
\newtheorem{cor}[thm]{Corollary}
\newtheorem{claim}[thm]{Claim}
\newtheorem {fact}[thm]{Fact}
\newtheorem*{thmstar}{Theorem}
\newtheorem{prop}[thm]{Proposition}
\newtheorem*{propstar}{Proposition}
\newtheorem {lem}[thm]{Lemma}
\newtheorem*{lemstar}{Lemma}
\newtheorem{conj}[thm]{Conjecture}
\newtheorem{question}[thm]{Question}
\newtheorem*{questar}{Question}
\newtheorem{ques}[thm]{Question}
\newtheorem*{conjstar}{Conjecture}
\theoremstyle{remark}
\newtheorem{rem}[thm]{Remark}
\newtheorem{np*}{Non-Proof}
\newtheorem*{remstar}{Remark}
\theoremstyle{definition}
\newtheorem{defn}[thm]{Definition}
\newtheorem*{defnstar}{Definition}
\newtheorem{exam}[thm]{Example}
\newtheorem*{examstar}{Example}

\newcommand{\pd}[2]{\frac{\partial #1}{\partial #2}}
\newcommand{\pp}{\partial }
\newcommand{\pdtwo}[2]{\frac{\partial^2 #1}{\partial #2^2}}
\newcommand{\od}[2]{\frac{d #1}{d #2}}
\def\Ind{\setbox0=\hbox{$x$}\kern\wd0\hbox to 0pt{\hss$\mid$\hss} \lower.9\ht0\hbox to 0pt{\hss$\smile$\hss}\kern\wd0}
\def\Notind{\setbox0=\hbox{$x$}\kern\wd0\hbox to 0pt{\mathchardef \nn=12854\hss$\nn$\kern1.4\wd0\hss}\hbox to 0pt{\hss$\mid$\hss}\lower.9\ht0 \hbox to 0pt{\hss$\smile$\hss}\kern\wd0}
\def\ind{\mathop{\mathpalette\Ind{}}}
\def\nind{\mathop{\mathpalette\Notind{}}}
\numberwithin{equation}{section}

\def\id{\operatorname{id}}
\def\Frac{\operatorname{Frac}}
\def\Const{\operatorname{Const}}
\def\spec{\operatorname{Spec}}
\def\span{\operatorname{span}}
\def\exc{\operatorname{Exc}}
\def\Div{\operatorname{Div}}
\def\cl{\operatorname{cl}}
\def\mer{\operatorname{mer}}
\def\trdeg{\operatorname{trdeg}}
\def\ord{\operatorname{ord}}

\newcommand{\m}{\mathbb }
\newcommand{\mc}{\mathcal }
\newcommand{\mf}{\mathfrak }
\newcommand{\is}{^{p^ {-\infty}}}
\newcommand{\QQ}{\mathbb Q}
\newcommand{\fh}{\mathfrak h}
\newcommand{\CC}{\mathbb C}
\newcommand{\RR}{\mathbb R}
\newcommand{\ZZ}{\mathbb Z}
\newcommand{\tp}{\operatorname{tp}}
\newcommand{\SL}{\operatorname{SL}}

\title{Model theory and machine learning}
\author{Hunter Chase}
\address{Department of Mathematics, UIC, Chicago IL}
\email{hchase2@uic.edu}

\author{James Freitag}
\address{Department of Mathematics, UIC, Chicago IL}
\email{freitagj@gmail.com}

\thanks{James Freitag was supported by NSF grant no. 1700095}

\maketitle

\begin{abstract} 
About 25 years ago, it came to light that a single combinatorial property determines both an important dividing line in model theory (NIP) and machine learning (PAC-learnability). The following years saw a fruitful exchange of ideas between PAC learning and the model theory of NIP structures. In this article, we point out a new and similar connection between model theory and machine learning, this time developing a correspondence between \emph{stability} and learnability in various settings of \emph{online learning.} In particular, this gives many new examples of mathematically interesting classes which are learnable in the online setting. 
\end{abstract}

\section{Introduction} 
The purpose of this note is to describe the connections between several notions of computational learning theory and model theory. The connection between \emph{probably approximately correct} (PAC) learning and the non-independence property (NIP) is well-known and was originally noticed by Laskowski \cite{laskowski1992vapnik}. In the ensuing years, there have been numerous interactions between the combinatorics associated with PAC learning and model theory in the NIP setting. Below, we provide a quick introduction to the PAC-learning setting as well as learning in general. Our main purpose, however, is to explain a \emph{new connection} between the model theory and machine learning. Roughly speaking, our manuscript is similar to \cite{laskowski1992vapnik}, but develops the connection between \emph{stability} and \emph{online learning}. 

That the combinatorial quantity of VC-dimension plays an essential role in isolating the main dividing line in both PAC-learning and perhaps the second most prominent dividing line in model-theoretic classification theory (NIP/IP) is a remarkable fact. This connection has been the subject of numerous works in recent years \cite{chernikov2013externally, guingonanip, johnson2010compression, livni2013honest}. In the setting of \emph{online learning} (described below), another combinatorial notion, the \emph{Littlestone dimension}, isolates the dividing line between learnability and non-learnability of a concept class. Given how well-studied the connection between model theory and the combinatorics associated with machine learning is, it is surprising that it hasn't been noticed until now that the same combinatorial quantity isolates what is perhaps the most prominent dividing line in classification theory (stable/unstable). 

Now we roughly describe the PAC setting, in part to contrast the setting with that of online learning. Given an infinite set $X$ with a probability measure $\mu$ on $X$ and a collection of measurable subsets of $X$, denoted by $\mc F$, one attempts to ``learn" a fixed but unknown $F \in \mc F$ by sampling from $X$. For some large $n$, $n$ elements of $X$ are randomly sampled, and the goal is to estimate the probability $\mu(F)$ by the proportion of elements of the sample which lie in $F$. For some $\epsilon >0$ fixed ahead of time, we say that the sample estimates the set $F$ \emph{$\epsilon$-well} if the proportion of elements of the sample which lie in $F$ is within $\epsilon$ of $\mu (F)$. The class $\mc F$ is \emph{learnable} if for any $\delta$ there is a large enough $n$ such that the measure of the samples of size $n$ (computed using the product measure $\mu^n$) which estimate the sample $\epsilon$-well is greater than $1-\delta$. Roughly, \emph{for large enough sample size, we can get arbitrarily high likelihood that a sample estimates the true probability arbitrarily well.} That is, for a large enough sample size, predictions are \emph{probably approximately correct.} It turns out that there is a purely combinatorial characterization of $\mc F$ being PAC-learnable (which remarkably does not depend on the distribution $\mu$); the collection $\mc F$ is PAC-learnable if and only if $\mc F$ has finite \emph{VC-dimension}. 

The connection to model theory is as follows: when $X$ is taken to be $\mc M$, a model of a first order theory $T$ and $\phi(x,y)$ is a formula in the language of $T$, we let $\mc F =\{\phi (\mc M, a) \, | \, a \in \mc M \}$. Then the VC-dimension of $\mc F$ is finite if and only if $\phi(x,y)$ is NIP. 

In the most straightforward (and restrictive) setup of online learning, we are given an infinite set $X$ (with no distribution) along with a collection $\mc F$ of subsets of $X$. The collection $\mc F$ is known to the learner. Fix some $F \in \mc F$ which is not known to the learner. Fixing some large $n$, there will be $n$ rounds. In round $i$, an element $x_i$ is selected, and the learner must predict the value of $1_F (x_t),$ that is, whether or not $x_i$ is in the unknown set $F$. We call the value of the learner's prediction $\hat y_i$. The goal of online learning is to minimize the number of mistakes made during these predictions $$\sum _{i=1} ^ n |\hat y_i - 1_F (x_i)|.$$ In this setting, there is no assumption about how the elements $\bar x = (x_1, \ldots , x_n)$ are chosen, and the choice of $x_{i+1}$ is allowed to depend on the predictions made by the learner in the previous rounds. One seeks to minimize the number of mistakes over all possible sequences of samples. This setting of computational learning often arises when the data becomes available in sequential order or the data is chosen by a process which is assumed to be adversarial to the learner (a process or opponent seeking to make the number of mistakes large). Variations on how the samples are chosen are possible as well; for instance, a certain limited amount of randomness is often injected into how the elements $x_i$ are chosen without moving the sampling back into the PAC context. 

It turns out that the number of mistakes that the best deterministic algorithm makes (over all possible samples) can be bounded in terms of a combinatorial quantity associated with the collection $\mc F$, the Littlestone dimension. When $X$ is taken to be $\mc M$, a model of a first order theory $T$, $\phi(x,y)$ is a formula in the language of $T$, and $\mc F =\{\phi (\mc M, a) \, | \, a \in \mc M \}$, the Littlestone dimension (also called thicket dimension) is precisely the Shelah 2-rank of $\phi (x,y)$, which is finite if and only if $\phi(x,y)$ is stable. A number of variants of this basic setup have much less restrictive assumptions (sometimes with a certain amount of randomness similar to the PAC setting) while also having the property that learnability is characterized by stability. In section \ref{OLSTAB} we will give an exposition of the various settings in which stability characterizes learnability. 

It seems surprising to the authors that the connection pointed out in the previous paragraph has not been previously noticed, but the following quote of \cite{rakhlin2011online} offers something of an explanation: 
\begin{quote} 
A reflection on the past two decades of research in learning theory reveals (in our somewhat biased view)
an interesting difference between Statistical Learning Theory and Online Learning. In the former, the focus
has been primarily on understanding complexity measures rather than algorithms... In contrast, Online Learning has been mainly
centered around algorithms.
\end{quote}
The dividing lines in model-theoretic classification theory are more naturally associated with combinatorial properties and the various complexity measures associated with PAC learning than with algorithms, and in the less restrictive online setups, the role of Littlestone dimension is perhaps somewhat more hidden than the role of VC-dimension in the PAC setup. 

The correspondence between online learnability and stability is similar to the correspondence between PAC learnability and NIP, but it should be mentioned that the fields (online learning and stability theory) are in rather different positions than in PAC learning correspondence with NIP. At this point, stability theory has been extensively developed, while at the time of \cite{laskowski1992vapnik}, the study of theories without the independence property was in its infancy, while PAC learning was much more developed. Various notions from PAC learning eventually played a big role in the development of structural results for NIP structures. In the case of the correspondence between stability and online learning, there seems to be more potential for the application of model theoretic ideas in online learning. For instance, in the final sentence of \cite{ben2009agnostic}, the authors mention that one of the main open questions in the theory is to close the gap between the lower bounds and upper bounds for the expected number of mistakes a learner makes in various online contexts, and that this question seems to have as a main obstacle a lack of interesting infinite concept classes with finite Littlestone dimension. Model theory offers a remedy for this obstacle; a great many mathematically interesting theories have been proven to be stable over the last forty plus years of classification theory, often with highly nontrivial proofs. So, following our discussion of online learning, we give some prominent examples of stable theories, giving various new examples of classes of finite Littlestone dimension. 

Now we describe the organization of this manuscript. In section \ref{general}, we describe the setting of computational learning in very general terms. In section \ref{PACNIP} we specialize to the PAC setting. In section \ref{OLSTAB}, we specialize to the setting of online learning before describing several variants. In the final section, we survey some stable theories, and use the connection pointed out earlier in the paper to give many new examples of classes with finite Littlestone dimension. 

\subsection{Acknowledgements} 
The authors would like to thank Siddharth Bhaskar, Alex Kruckman, Dimitrios Diochnos, Dave Marker, Lev Reyzin, Dhruv Mubayi, Maryanthe Malliaris, and Gyorgy Turan for useful suggestions and conversations during the preparation of this article. 

\section{Machine learning generalities} \label{general} 
In this section, we describe the generalities of machine learning, in quite a general setup, while mentioning the cases of particular interest to us. Let $Y$ be a set, which we will call the set of \emph{labels}. Let $Y'$ be another set, which we will refer to as the \emph{predictions.} Fix a function $$L: Y \times Y' \rightarrow \m R _{\geq 0}$$ which we call the \emph{loss function}. 

\begin{rem} 
The most common setup occurs when $Y=Y'=\{0,1\}$ and $L(y,y')= |y-y'|.$ 

Another common example occurs when $Y=Y'=I \subseteq \m R,$ with $I$ a bounded interval. In this case, a common loss function is given by $L(y,y')= (y-y')^2.$ Settings in which $Y,Y' \subset \m R$ are sometimes called margin-based. These settings are less natural to connect directly to model theory, though it might make sense to study margin-based machine learning in the context of continuous model theory \cite{yaacov2006model}. 
\end{rem}

Let $X$ be another set, which we call the set of \emph{examples} (also sometimes called \emph{inputs} or \emph{instances}). A \emph{concept} is a map $c: X \rightarrow Y$. In the example given above with $Y= \{0,1\}$, a concept is simply a subset of $X$. A concept class $\mc C$ is a collection of concepts. 

Fix some concept $c$. The learner will make a series of predictions about a sample of inputs from $X$ by selecting a prediction $\hat y_i$ for the label of each element $x_i$ from the sample. The learner incurs a loss for each element $x_i$ of the sample, by evaluating $L(c(x_i), \hat y_i)$. If the elements of the sample are indexed by the set $I$, then the total loss incurred is given by $$\sum _{i \in I} L(c(x_i), \hat y_i). $$

The goal of the learner is always the same---minimize the total loss coming from making predictions about a series of elements of $X$. Besides the objects described above, the differences in various settings of learning theory are derived from the assumptions about what data the learner has available and how the elements of the sample are chosen.

\section{PAC-learning and NIP} \label{PACNIP}
In this section, we will quickly explain the connection between PAC learning and NIP. Our presentation essentially follows \cite{guingonanip}. Fix a concept class $\mc C$ on a set $X$ with $Y=Y' = \{0,1\}$. Let $\mc C_{fin} = \{ f |_ Y \, | \, Y \subset X , \, Y \text{ finite, } f \in \mc C\}$. Let $\mu$ be a probability measure on $X$ such that each element of $\mc C$ is measurable. We will think of the learner as having complete knowledge of the elements of $\mc C$, and the elements for a sample being drawn randomly with respect to the distribution given by $\mu$. 

Let $G:\mc C_{fin} \rightarrow 2^X$ be a function. Let $\bar a = (a_1, \ldots , a_n)$. Define $$err_{\mu} (G, f, \bar a):= \mu ( \{ c \in X \, | \, f(c) \neq G( f|_{\bar a}) (c) \} ).$$ Here one should think that $G$ is a function being used to generate predictions, while the error is the probability that the next prediction is incorrect. 

We say that $\mc C$ is \emph{probably approximately correct learnable} (PAC-learnable) if there is a $G: \mc C_{fin} \rightarrow 2^X$ such that for all $\epsilon >0$ and all $\delta >0$, there is $N_{\epsilon, \delta } \in \m N$ such that for all $f \in \mc C$, and all $\mu $ on $X$ such that all elements of $\mc C$ measurable, $$\mu^{N_{\epsilon, \delta}} \left( \{ \bar a \in X^{N_{\epsilon, \delta }}  \, | \, err_\mu (G, f, \bar a) > \epsilon \} \right) < \delta ,$$ where $\mu^{N_{\epsilon, \delta}}$ is the product measure. 
That is, the probability that the error is high (bigger than $\epsilon$) is small (less than $\delta$). Supposing that the class $\mc C$ is PAC-learnable, there is a minimal $N_{\epsilon, \delta}$ for which the inequality holds, which is called the \emph{sample complexity}. 

The following theorem establishes the connection between VC-dimension and PAC-learnability: 

\begin{thm} Let $\mc C $ be a concept class on $X$. Then the following are equivalent: 
\begin{enumerate} 
\item $\mc C$ has finite VC-dimension. 
\item $\mc C$ is PAC-learnable, and $$N_{\epsilon, \delta} \leq  max \left\{ \frac{4}{\epsilon} \log _2 \left(\frac{2}{\delta} \right) , \frac{8d}{\epsilon} \log _2 \left( \frac{13}{\epsilon} \right) \right\} .$$
\end{enumerate} 
\end{thm}

In fact, even more is true---if $\mc C$ is PAC-learnable with sample complexity $N_{\epsilon, \delta}$, then one can show that the expected value of the function $\bar a \mapsto err_\mu (G, f , \bar a)$ is bounded by $ \delta + \epsilon (1- \delta).$ 

In the years since Laskowski's paper \cite{laskowski1992vapnik}, connections between the VC theory and NIP have developed extensively with important notions from VC-theory adapted to the model-theoretic setting and vice versa \cite{chernikov2013externally, guingonanip, johnson2010compression, livni2013honest}.

\section{Online learning and stability} \label{OLSTAB} 
The initial setting of online learning which we describe is due to Littlestone \cite{littlestone1988learning}; the particular setting received relatively little attention, perhaps due to the very strong assumptions (\cite{littlestone1988learning} is in fact famous for several other contributions). Littlestone's work was generalized in various ways in the ensuing years, with the assumptions being significantly weakened. We will begin with the original setup of \cite{littlestone1988learning}, and eventually describe two settings laid out in \cite{ben2009agnostic}. First, we set up some of the combinatorial notions pertinent in each of the settings we consider.  

The next several definitions follow the notation and terminology of Bhaskar \cite{bhaskar2017thicket}. 

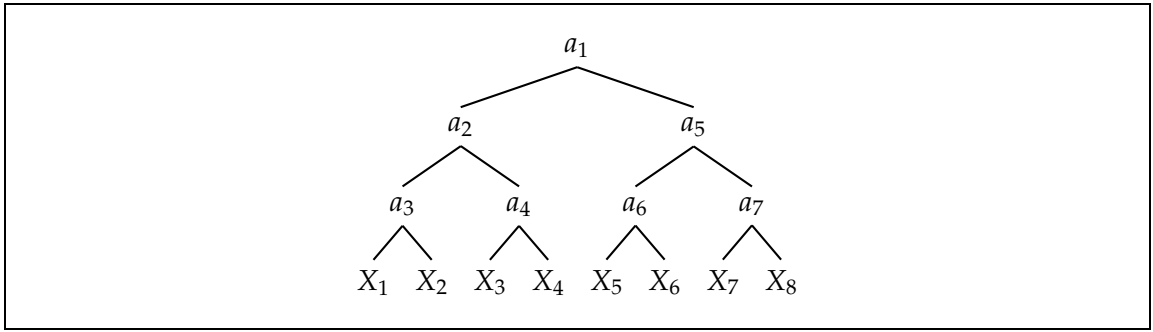
\begin{figure}[h]
\centering
\begin{framed}
\begin{tikzpicture}
\tikzset{edge from parent/.style={draw,thick}}
\Tree [.$a_1$ [.$a_2$ [.$a_3$ [.$X_1$ ] [.$X_2$ ] ] [.$a_4$ [.$X_3$ ] [.$X_4$ ] ] ]
[.$a_5$ [.$a_6$ [.$X_5$ ] [.$X_6$ ] ] [.$a_7$ [.$X_7$ ] [.$X_8$ ] ]]]
\end{tikzpicture}
\end{framed}
\caption{A binary element tree of height three. Here $a_ i \in X$ and $X_i \in \mc F$. The leaf labeled with $X_4$ is well-labeled if and only if $a_1 \in X_4$ and $a_2, a_4 \notin X_4$. For all other $a_i$, there is no requirement about membership in $X_4$.} \label{BDT}
\end{figure}

\begin{defn} A \emph{binary element tree of height $h$}, denoted by $\mc T _h$, is a rooted complete binary tree of height $h$ whose non-leaf vertices are labeled by elements of the set $X$ and whose leaves are labeled by elements of $\mc C$ (see Figure \ref{BDT}).
\end{defn} 
For the following definitions, fix a binary element tree of height $h$. 
\begin{defn} 
A vertex $v_1$ is \emph{below} a vertex $v_2$ if $v_2$ lies on the (unique) path from $v_1$ to the root of the tree. We say that $v_1$ is \emph{left-below} $v_2$ if $v_1$ is below $v_2$ and the first edge along the path from $v_2$ to $v_1$ goes down and to the left. The notion of \emph{right-below} is defined analogously. When a vertex labeled by $b$ is left-below a vertex labeled by $a$, we write $a <_L b$. Similarly, when a vertex labeled by $b$ is right-below a vertex labeled by $a$, we write $a <_R b$.
\end{defn}

\begin{defn}
A leaf, labeled by $Y \in \mc C$ is said to be \emph{well-labeled} if for each vertex above $Y$, say labeled by $a$, 
$$a \in Y \text{ if and only if }  a <_L Y.$$  
\end{defn} 

\begin{defn}  The \emph{thicket shatter function}  $\rho_{\mc F}:\ZZ^{\ge 0} \rightarrow \ZZ^{\ge 0}$ is defined by letting $\rho_{\mc F}(n)$ be the maximum number of well-labeled leaves on a binary element tree of height $n$, $\mc T_n$, whose leaves are labeled with elements of $\mc F$. The \emph{thicket dimension} $Ldim(\mc F)$ is the maximum integer $n$ such that $\rho _\mc F (n) = 2^n$, or else $Ldim(\mc F)= \infty$. 
\end{defn} 

Thicket dimension has appeared in at least several other contexts under different names; in fact Bhaskar \cite{bhaskar2017thicket} was aware of the terminology and definitions of \cite{Shelah}, which we reproduce next: 

\begin{defn} Let $\mc M$ be a monster model of a complete $\mc L$-theory. Fix a consistent partial type $ \pi (x)$ and a partitioned formula $\phi (x;y)$. Then the ordinal $R(\pi, \phi , 2)$, called the Shelah 2-rank, is defined as follows: 
\begin{itemize} 
\item $R(\pi, \phi , 2) \geq 0$. 
\item For any limit ordinal $\lambda,$ $R(\pi, \phi ,2) \geq \lambda $ if $R(\pi, \phi , 2) \geq \alpha $ for all $\alpha < \lambda$. 
\item For any ordinal $\alpha$, $R(\pi, \phi ,2) \geq \alpha+1 $ if there is some $\phi (x,a)$ such that $R(\pi \cup \{\phi (x,a) \} , \phi, 2) \geq \alpha$ and $R(\pi \cup \{\neg \phi (x,a) \} , \phi, 2) \geq \alpha$.
\end{itemize} 
\end{defn} 

In general, $R(\pi,\Delta,2)$ can also be defined for a finite collection of formulas $\Delta$, but this case can be shown to reduce to the case of a single formula.
The formula $\phi(x,y)$ is \emph{stable} if and only if $R(\emptyset, \phi ,2)$ is finite \cite{Shelah}; a theory is stable if every formula is stable. It is reasonably clear that the $R(\pi,\phi,2)$ is the thicket dimension of the set system on $\mc M^{|y|}$ given by the collection of sets $\{ \phi(b,\mc M) \, | \, b \in \pi (\mc M)\}$; for more details, see \cite{bhaskar2017thicket}. 

The thicket dimension also appears for the first time in the context of learning theory in \cite{littlestone1988learning}; the quantity came to be called the \emph{Littlestone dimension} \cite{ben2009agnostic}.

\subsection{The realizable case} Fix a set system $\mc C $ on a set $X$. 
Assume that $Y=Y' = \{0,1\}$ and the loss function for a prediction $\hat y$ and concept (that is, a set) $X$ on input $x$ is given by $|\hat y - 1_X(x)|.$ Over all possible algorithms, we seek to minimize our loss, that is, the number of mistakes we make over $n$ rounds of predictions. In the \emph{realizable} case, we assume that $X \in \mc C$, so that the true concept is among the set of concepts $\mathcal{C}$ accessible to the learner. There are \emph{no} assumptions on the choices of the instances $x_t$. The goal is to minimize the worst case number of mistakes made by our predictions over all possible samples of the instances and choice of the concept. So, we seek to bound $$M = \max _{X \in \mc C} \max _{\bar x = (x_1, \ldots , x_n) } \sum _{t=1}^n | \hat y_t - 1_X (x_t)|,$$ where $\hat y_t$ is chosen by some deterministic algorithm. 

For applications and purposes of discussing the bounds, one often views the entity selecting the instances $\bar x$ as antagonistic to the learner---and in our current simplified setting, bounding the worst case number of mistakes bounds the actual number of mistakes made when the antagonistic sampling entity has perfect information about the prediction process.

\begin{thm} \cite{littlestone1988learning} 
The worst case number of mistakes of any deterministic algorithm in the online learning setting with concept class $\mc C$ is at least the Littlestone dimension of $\mc C$, and there is an algorithm that makes at most this many mistakes. 
\end{thm}

\begin{rem} The algorithm which minimizes the number of worst-case mistakes in the above setting is referred to as the Standard Optimal Algorithm (SOA), and we describe it briefly here. Begin with $V_0 = \mc C$. At each stage, the learner inductively defines $V_i$. At stage $t$, the learner receives $x_t$, and sets, for $r=0,1$, $$V_t ^{(r)} := \{ X \in V_ {t-1} \, | \, 1_X ( x_t) =r \}.$$ The learner predicts $\hat y_t=r$ which maximizes the Littlestone dimension of $V_t^{(r)}$ (ties are predicted in some fixed manner, say $\hat y_t=0$ in the case of a tie). Then the learner gets the value of $1_X(x_t)$ and realizes whether a mistake has been made. At this point, set $V_t = V_t^{1_X(x_t)}.$ 

The essential point here is that if a mistake is made, it must be the case that the Littlestone dimension of $V_t$ is strictly less than the Littlestone dimension of $V_{t-1}$ (proving this is an easy exercise). Of course, this bounds the total number of mistakes which the algorithm can ever make under any choice of $\bar x$ by the Littlestone dimension. 
\end{rem}

\subsection{Learning from experts}
The case in which we assume that the learner has access to true concept $X \in \mc C$ is often referred to as the \emph{realizable} case of online learning. For various applications, this assumption is too strong (as are other assumptions from the previous subsection which we will deal with in later sections). In this section, we will explain a context of online learning which removes the realizability assumption.

The goal again is to minimize mistakes, but here, the minimization will be relative to a particular class of $\{0,1 \}$-valued functions, which we will call $\mc H .$ That is, we wish to minimize, for any sampling of instances, $\bar x = (x_1, \ldots , x_T)$, the difference between the number of mistakes made by the learner and the minimal number of mistakes made by any of the functions in $\mc H.$ So, in this case, the loss function is taken to be 
$$\sum |\hat y_t - y_t| - \min_{h \in \mc H} \sum |h(x_t) - y_t|.$$ 
Here one often thinks intuitively that the functions in $\mc H$ are experts making predictions, and the learner's job is to choose which expert's prediction to believe.

Littlestone and Warmuth \cite{littlestone1994weighted} consider this problem in the case that $\mc H$ is finite via a probabilistic weighted majority algorithm. We will now describe their algorithm. At the outset, each of the $N$ many experts $\{f_i \}_{i=1}^N=\mc H$ is assigned weight $1$, and the weight of expert $i$ at stage $t$ will be denoted by $w_i^{t}$. We fix the learning rate $\eta > 0$, which dictates how much we discount the weight of an expert for providing incorrect advice. At each stage, the learner receives the expert advice, $(f_1(x_t), \ldots , f_N (x_t))$, a tuple in $\{0,1\}^N$. The learner predicts $1$ with probability $$p_t = \frac{1}{\sum_{i=1}^N w_i^{t-1}} \sum_{i=1}^N w_i ^{t-1} f_i (x_t).$$ Then once the actual value $y_t$ is revealed, the weights are updated via: $w_i^t = w_i^{t-1} e^{- \eta \cdot | f_i (x_t) -y_t|}.$ That is, those experts who were wrong see their weight drop by a factor of $e^{-\eta}$. 

The expected value of the loss function of their algorithm with a sample of size $T$ is 
$$\sum_{t=1}^T E(|\hat y_t -y_t|)-\min_{h \in \mc H} \sum_{t=1}^T |h(x_t)-y_t| \leq \sqrt{\frac{1}{2} \ln(N)T}.$$ 
Here, the assumption that $\mc H$ is finite is often too strong for applications, however, \cite{ben2009agnostic} generalize the setup to the case in which $\mc H$ is infinite, but of finite Littlestone dimension, proving: 

\begin{thm} 
There is an algorithm such that for all $h \in \mc H$ and any sequence of instances $\bar x = (x_1, \ldots, x_T),$ 
$$\sum_{t=1}^T E(|\hat y_t -y_t|)-\min_{h \in \mc H} \sum_{t=1}^T |h(x_t)-y_t| \leq \sqrt{\frac{1}{2} Ldim (\mc H ) \cdot T \ln (T) }.$$
\end{thm}

In \cite{ben2009agnostic} it is also shown that no algorithm (even allowing randomization) can achieve an expected bound better than $\sqrt{\frac{1}{8} Ldim(\mc H) T}.$ Closing the gap between the lower and upper bounds for the loss function (sometimes called regret in this context) is one of the main open problems mentioned in \cite{ben2009agnostic}, where the authors remark that there are few known interesting examples of infinite classes with finite Littlestone dimension. 

\subsection{Bounded stochastic noise}
Suppose that we work in the general setup from the previous section (again, not assuming realizability), but with a difference in the way we generate labels and measure mistakes. Suppose that there is a function $h \in \mc H$ such that the labels $y_1, \ldots , y_ T$ are independent $\{0,1 \}$-valued random variables with the property that for all $t$, $Pr(h (x_t ) \neq y_t ) \leq \gamma$ with $\gamma \in (0,\frac{1}{2})$. This value $\gamma$ will be called the noise rate. 

In this setting, one seeks to minimize the difference between the predictions and the output of the noisy function on the samples:
$$E \left( \sum _{t=1}^T  | \hat y_t - y_t | \right).$$ 
Note here that there are two sources of randomness---the choices of the algorithm may be randomized and the labels $y_t$ are random variables. The expectation is taken with respect to both of these. 

\begin{thm} For any concept class $\mc H$, and any $\gamma \in [0,\frac{1}{2})$, there is an algorithm (possibly randomized) so that for any $h \in \mc H$, and a sequence of examples $(x_1, y_1), \ldots , (x_T, y_T)$ with each $y_t$ a random variable as described above, 
$$E \left( \sum _{t=1}^T | \hat y_t - h (x_t)| \right) \leq \frac{Ldim(\mc H) \cdot \ln (T) }{1- 2 \sqrt{\gamma (1- \gamma)} }.$$ 
\end{thm} 
That is, the expected number of mistakes grows only logarithmically in the sample size. In \cite{ben2009agnostic}, the authors give an example of a class $\mc H$ which shows that the left hand side of the inequality in the theorem is bounded below by $\Omega (Ldim(\mc H) \cdot \ln (T) ).$

\section{Stability theory} 
In this section, we use stability theory to point out various mathematically interesting examples of classes which have finite Littlestone dimension. We will assume some basic familiarity with first order logic, but we provide some reminders for the non-model theorist for whom this section is written.

Fix some complete theory $T$ in a language $\mc L$ and let $\mc M$ be a monster model of $T$. The non-model theorist can simply loosely assume that $\mc M$ is a very large structure in which over a small subset $A$ (say of cardinality at most $\kappa$) for any tuple $c $ in any model of $T$ containing $A$, there is some $b \in \mc M$ such that $\tp (c/A) = \tp (b/A)$. Here $\tp (c/A)$ denotes the collection of all first order formulas in the language $\mc L$ with parameters from $A$ which are satisfied by $c$.

For $n \in \m N$, the space of types of $n$-tuples of $\mc M$ over some subset $A \subset \mc M$ is denoted by $S_n (A)$. It comes naturally equipped with a topology in which the basic open sets correspond to first order formulas with parameters in $A$. Rather than considering all formulas, sometimes it is natural to restrict to the $\phi$-type of a tuple, denoted $\tp_\phi (c /A)$, the collection of instances of $\phi$ with parameters in $A$ which hold of $c$.
When $\phi (x;y)$ is a formula, the space of $\phi$-types over $A$ (treating the variables $y$ as parameters) is denoted by $S_\phi (A)$. 

The theory $T$ is called $\kappa$-stable if for every set $A \subseteq \mc M$ with $|A| \leq \kappa$, we have $|S_n (A) | \leq \kappa$ for all $n \in \m N$. The theory is \emph{stable} if it is $\kappa$-stable for some $\kappa \geq |T|$. Part of the utility of the notion is that it can be characterized in several disparate ways (this is not an exhaustive list): 
\begin{fact} \cite{Shelah} The following conditions are equivalent: 
\begin{enumerate} 
\item $T$ is $\kappa$-stable for some $\kappa$.
\item For any countable set $A \subset \mc M$, $S_{ \phi } (A)$ is countable. 
\item Every formula $\phi (x;y)$ has finite Shelah $2$-rank---that is, $R( \emptyset , \phi, 2)$ is a finite ordinal (recall that Shelah $2$-rank is equal to Littlestone dimension). 
\item No formula $\phi (x;y)$ has the \emph{order property}. A formula $\phi(x;y)$ has the order property if there are tuples $(a_1, b_1), (a_2, b_2), \ldots $ from $\mc M$ so that $\mc M \models \phi (a_i ;b_j)$ if and only if $i \leq j$. 
\end{enumerate} 
\end{fact} 

When $\kappa$ in the first condition of the above definition is be taken to be $\aleph _0$, the theory is (somewhat enigmatically) called $\omega$-stable. Not every stable theory is $\omega$-stable, even when making strong assumptions about various aspects of the language or structure. For instance, the theory of the integers where the language consists of the additive group operation as a binary function is stable, but not $\omega$-stable. 

Stability is one of the dividing lines (probably the most prominent one) which in certain contexts, model-theorists view as the border between ``tame" and ``wild" structures; stability allows for the development of various structural results, which are (often provably) impossible in the case of unstable theories. Stability has various non-obvious interactions with algebraic structure, and understanding these interactions has been the subject of a huge amount of model theoretic work over the past fifty years (for instance, there is a deep structure theory of stable groups \cite{Poizat}). 
 
Consider the concept class $\mc C_\phi$ on $\mc M^{|y|}$ given by the collection of sets $\{ \phi(b,\mc M) \, | \, b \in \pi (\mc M)\}$. The theory $T$ is stable precisely if each concept class of this form has finite Littlestone dimension (see section \ref{OLSTAB} for an explanation). 

We will elaborate on condition (4). Given a class $\mc C _ \phi $, there is a natural bipartite graph $G_\phi$ associated with any concept class. The sets of vertices consist of 1) the elements of the underlying set and 2) concepts, with an edge between an element and a concept if and only if the element is in the concept. Finite Littlestone dimension of the concept class $\mc C_ \phi$ is equivalent to there being an upper bound on the size of any half-graph which appears as an induced subgraph of $G_ \phi$. 

\subsection{Examples of notable stable theories} 
We now make a list (very far from comprehensive) of some notable stable theories and offer some explanation of the set systems (families of definable sets) which arise in the various settings. From our list, many mathematically interesting classes $\mc C_ \phi$ with finite Littlestone dimension can be obtained. 
\begin{enumerate} 
\item $ACF$, the theory of algebraically closed fields. By quantifier elimination for algebraically closed fields, the concept classes which appear as $\mc C_\phi$ in the theory of algebraically closed fields are precisely the uniform families of affine constructible sets. That is, when $f: V \rightarrow W$ is a rational map (everything defined over some fixed algebraically closed field), the corresponding family of constructible sets is the collection of fibers of the function $f$. More concretely, one can think of such a family as being given by solutions sets of families of polynomial equations and inequations:
$$f_1(x, a)  = f_2 (x,a) , \ldots , f_n(x,a)  =0, f(x,a) \neq 0$$ where $x$ is a tuple of indeterminates and $a$ is a tuple which varies over some constructible subset of $\m A^{|a|}$. 

\item $DCF_0$, the theory of differentially closed fields of characteristic zero, was first investigated by Robinson \cite{Robinson}, and Blum \cite{Blum} gave an elegant axiomatization from which it was straightforward to notice that the theory is stable. See \cite{MMP} for a more comprehensive discussion of $DCF_0$, as we will be brief here. Differentially closed fields are universal domains for algebraic differential equations; that is, if a system of equations has a solution in some field of functions, it already has a solution in the differential closure of the field generated by the coefficients of the equations. By quantifier elimination for differentially closed fields, the concept classes which appear as $\mc C_\phi$ in the theory of differentially closed fields are precisely the uniform families of constructible sets in the Kolchin topology (boolean combinations of the zero sets of algebraic differential equations). That is, when $f: V \rightarrow W$ is a differential rational map between affine constructible sets $V,W$ in the Kolchin topology (everything defined over some fixed differentially closed field), the corresponding family of constructible sets is the collection of fibers of the function $f$. Such a family is alternatively given by a collection of differential equations and inequations
$$f_1(x, a)  = f_2 (x,a) , \ldots , f_n(x,a)  =0, f(x,a) \neq 0$$ where $x$ is a tuple of indeterminates from $\mc M \models DCF_0$ and $a \in \mc M$ is a tuple which varies over some Kolchin-constructible subset of $\m A^{|a|}$.

\item The theory of separably closed fields with characteristic $p \neq 0$ and fixed degree of imperfection $e \in \m N$ (which we will describe here) is complete and was shown to be stable by Wood \cite{wood1979notes}. When a field $F$ of characteristic $p$ is closed under separable extensions, we say $F$ is separably closed. A set $B \subseteq F$ is a \emph{$p$-basis} of $F$ if the collection of products of powers of elements of $B$ of degree at most $p-1$ forms a basis for $F$ as an $F^p$-vector space. The cardinality of such a set $B$ is called the degree of imperfection of $F$ (which we assume to be finite). Now let $\{a_1, \ldots , a_e \}$ be a $p$-basis of $F$, and let $\{m_1, \ldots , m_{p^e} \}$ be the collection of monomials in $\{a_1, \ldots , a_e \}$ of degree at most $p-1$ in each element. Every element of $F$ can be written uniquely in the form $$x = \sum _{i=1} ^ {p^e} x_{(i)}^p m_i$$ where $x_i \in F$. For each element $x_i$ in the above sum, we can repeat the process, writing $$x_{(i)} = \sum _{j=1} ^ {p^e} x_{(i,j)}^p m_i.$$ Naturally, one can continue to iterate this process, defining $x_{\sigma}$ for any $\sigma $ a finite tuple of elements from $\{1, \ldots , p^e \}$. Let $\lambda _\sigma$ be the unary function $x \mapsto x_{\sigma}$. 

Let $\mc L_{p,e} $ be the language $\{+,-,\cdot , ^{-1}, 0,1\} \cup \{a_1, \ldots , a_e\} \cup \{ \lambda _ \sigma \, : \, \sigma \in (p^e )^{<\omega} \}.$ The theory of separably closed fields of characteristic $p$ with degree of imperfection $e$ eliminates quantifiers in the language $\mc L_{p,e}$. So, in one variable, definable sets correspond to boolean combinations of the zero sets of ideals in $F[x, \lambda _\sigma (x) ]_{\sigma \in (p^e)^ {\leq n}},$ for some $n$. 

\item Let $X$ be a compact complex manifold. Consider the structure $\mc A (X)$ where the basic relations are the complex analytic subsets of $X^n$ for any $n \in \m N$; we call a subset $A \subseteq X^n$ complex analytic if it is, for any point $p \in X^n$ there is a neighborhood $U$ of $p$ such that $A \cap U$ is given by the zero set of some fixed finite number of holomorphic functions on $U$. The model theory of compact complex manifolds began with Zilber's observation \cite{ZilberCCM} that if one adds as a relation all complex analytic subsets of $X^n$ for all $n$, then the induced structure is stable. For an overview of the model theory of compact complex manifolds, see \cite{moosa2010model}.

\item Let $R$ be a ring and $\mc L_R$ be the language of right $R$-modules, consisting of a symbol for addition and a unary function $f_r$ for each $r \in R$, which is interpreted as scalar multiplication by $r$. Let $T$ be any complete theory of right $R$-modules in the language $\mc L_R$. By a result of Baur \cite{baur1976elimination}, every formula $\phi( x)$ is equivalent to a boolean combination of positive primitive formulas, that is, formulas of the form $\exists y \psi (x,y) $, where $\psi $ is a conjunction of atomic formulas. In particular, every definable subset of an $R$-module $M$ is a boolean combination of cosets of positive primitive definable subgroups of $M$. An abelian group can be viewed as a $\m Z$-module, and from this characterization of definable sets, it is not hard to show that every abelian group has a stable theory in the language of groups.

\item The theory of the nonabelian free group $T_{fg}$ in the language of groups was shown to be stable by Sela \cite{sela2006diophantine} (Sela shows the same for any torsion-free hyperbolic group). Every formula in the language of groups is, modulo the theory of the free group, equivalent to a $\forall \exists$-formula. The strategy of the proof is complicated and is developed by Sela over a series of seven previous papers; see \cite{sela2006diophantine} for complete references. 

\end{enumerate}

\bibliography{Research}{}
\bibliographystyle{plain}
\end{document}